%% file: ad2.tex
\documentclass{article}
\usepackage{amssymb,amscd,amsmath, amsthm, epsfig}

\def\Arr{\mathop{\rm Arr}}
\def\AArr{\mathop{\rm AArr}}

\def\MacP{\mathop{\rm MacP}}

\def\cross{\mathop{\rm cr}}

\def\Z2{{\mathbb Z}_2}              %integers mod 2
\def\Rn{{\mathbb R}^n}
\def\Rk{\mathbb R^k}
\def\Rinfty{{\mathbb R}^\infty}

\def\R{\mathbb R}
\def\Nn{\mathbb N}

\def\Mm{{\cal M}}

\newtheorem{theorem}{Theorem}[section]
\newtheorem{lemma}[theorem]{Lemma}

\newtheorem*{main}{Main Theorem}

\theoremstyle{definition}
\newtheorem{defn}[theorem]{Definition}

\def\beginwval#1#2{\bgroup% change to \the#1 should be local
  \edef\@savecount{\the\value{#1}}% save old value of counter
  \expandafter\def\csname the#1\endcsname{#2}% new value
  \begin{#1}}% begin environment
\def\endwval#1{\setcounter{#1}{\@savecount}\end{#1}\egroup}

\begin{document}
\title{There is no tame triangulation of the infinite real Grassmannian}
\author{Laura Anderson\footnote{Partially supported by grants from the
National Science Foundation.} \and James F. Davis$^*$} 
\date{}

\maketitle

\begin{abstract}
 We show that there is no triangulation of the infinite
real Grassmannian $G(k,\Rinfty)$ nicely situated with respect to 
the coordinate axes. In terms of matroid theory, this says there
is no triangulation of $G(k,\Rinfty)$  subdividing the matroid
stratification.  This is proved by an argument in projective
geometry, considering a specific sequence of 
configurations of points in the plane.
\end{abstract}
 
The Grassmannian $G(k,\Rn)$ of $k$-planes in $\Rn$ is a smooth manifold, hence can be triangulated.  Identify $\R^n$ as a
subspace of $\R^{n+1}$, and let $\Rinfty$ be the union (colimit) of the $\Rn$'s.  The Grassmannian $G(k,\Rinfty)$ is infinite
dimensional; it is unclear whether it can be triangulated for $k \geq 3$.  We are interested in triangulations which are nicely
situated with respect to the coordinates axes. Such triangulations are of interest in combinatorics in the context 
of matroid theory; see Section~\ref{context}.

\begin{defn}  A triangulation of $G(k,\Rn)$ or  $G(k,\Rinfty)$ is \textbf{tame} if for every simplex $\sigma$, for every pair
of $k$-planes $V,W \in \text{int } \sigma$, and for every vector $v \in V$, there is a vector $w \in W$ so that for all of the
standard basis vectors $e_i$,
$$v \cdot e_i = 0 \Leftrightarrow w \cdot e_i = 0
$$
\end{defn}

Using triangulation theorems from real algebraic geometry, it is not difficult to prove the following theorem (see \cite{Mac},
\cite{Z2}).

\begin{theorem}  For every $k$ and $n$, there is a tame triangulation  $T_{n,k}$.  Furthermore for $n' \leq n$ and $k'\leq k$,
the triangulation $T_{n,k}$ restricts to a subdivision of $T_{n',k'}$
\end{theorem}

This theorem does not lead to a triangulation of $G(k,\Rinfty)$, because perhaps one would have to infinitely subdivide
$G(k,\Rn) \subset G(k,\Rinfty)$.  

\begin{main}  There is no tame triangulation of $G(3,\Rinfty)$.
\end{main}

It follows immediately that there is no tame triangulation of $G(k,\Rinfty)$ for $k \geq 3$.

We first rephrase the main theorem in terms of matroids and oriented 
matroids and give some very basic context.
Section \ref{main} of the paper gives the proof of the main theorem, Section~\ref{generalize} gives generalizations of the
main theorem to more general  subdivisions and more general stratifications. Section~\ref{context}
discusses matroid bundles and the MacPhersonian, which gives the context for this paper.  The last three sections can be read
independently of one another.

\section{Matroid stratifications}\label{matroid}  

The motivation for tameness comes from matroids and oriented matroids, which are combinatorial abstractions of linear
algebra.  We don't give the definition here (see \cite{Z2} for the definition, and \cite{BLSWZ} for the full story), but
simply state that an oriented matroid on a set
$E$ is a subset of all functions from $E$ to the three-element set $\{+,-,0\}$,  where the subset satisfies certain axioms. A
similar definition for a matroid can be given as a subset of all functions from $E$ to $\{1,0\}$.  Any oriented matroid determines a matroid by
identifying $+$ and $-$ with 1. The functions in the
(oriented) matroid are called covectors.  An (oriented) matroid has a rank associated to it, and the MacPhersonian
$\MacP(k,n)$ is the set of all rank $k$ oriented matroids on $E = \{1,2, \dots, n\}$.  (The MacPhersonian has a natural
partial order, but this does not play a role in the proof of our main theorem.)   
Let $\{e_1,e_2,\ldots,e_n\}$ denote the unit coordinate vectors in $\Rn$.
For any $V\in G(k,\Rn)$ there is an associated rank $k$ oriented matroid
on $E=\{1,2,\ldots, n\}$ whose set of covectors is the set of all
$$\{i \mapsto \text{sign}\langle v, e_i \rangle\}_{v\in V}.
$$
There is also a matroid associated to the
arrangement of vectors, where $i \mapsto 1$ if and only if $\langle v, v_i \rangle $ is non-zero. 

\begin{defn}  The \textbf{(oriented) matroid stratification of $G(k,\Rn)$} is the partition in which $V,W$ are in the same
element of the partition if and only if they determine the same (oriented) matroid.  Elements of the partition are called
(oriented) matroid strata.
\end{defn}
 
We now see that a tame triangulation is nothing more than a triangulation which refines the matroid stratification, that is,
every stratum is a union of interiors of simplices. 
If $M_1$ and $M_2$ are distinct rank $k$ oriented matroids on $\{1,2,\ldots,n\}$ with the same underlying matroid, then the strata of $M_1$ and $M_2$ 
lie in disjoint open subsets of $G(k,\Rn)$. Thus any connected subset of a 
matroid stratum is contained in an oriented matroid stratum, and so any tame triangulation must refine the oriented matroid stratification. Our main proof will show that no triangulation of $G(3,\Rinfty)$ refines the oriented matroid 
stratification.

The (oriented) matroid stratification is interesting geometrically and combinatorially, and has been studied extensively
\cite{BLSWZ}.  We note two results.  First is the observation of Gelfand-Goresky-MacPherson-Serganova that the matroid
stratification is precisely the coarsest common refinement of all of the Schubert cell decompositions given by the standard
basis and permutations of the standard basis.  The second result is the theorem of Mn\"ev \cite{Mn} that the oriented matroid
strata can have arbitrarily ugly homotopy type, i.e., for any semialgebraic set $S$, there is an $n$ and an oriented
matroid stratum of
$G(3,\Rn)$ for some $n$ having the homotopy type of $S$.

 Our interest in this question
arose in considering the theory of \emph{matroid bundles}
(cf.~\cite{bundles}), a combinatorial model for real vector bundles. 
The relationship between $G(k,\Rinfty)$ and $\MacP(k,\infty)$
 is a critical question in the theory. The results in this paper arose as
a revelation that extending arguments from the finite Grassmannians
to the infinite Grassmannian is
harder than one might anticipate.

\section{Proof of the Main Theorem}\label{main}

The proof involves constructing a  sequence of specific elements of $G(3,\Rinfty)$ which, 
assuming a tame triangulation, leads to an infinite number of simplices in a compact space $G(3,\R^8)$, and hence to a 
contradiction.  

However, it is rather difficult to visualize $3$-planes in $\Rn$, so instead we work
with arrangements of vectors in $\R^3$.   Let 
$\Arr(k,n)$ be the space of all spanning $n$-tuples
$(v_1, \ldots, v_n)$ of vectors in $\Rk$.  Include $\Arr(k,n) \subset \Arr(k,n+1)$ by adding the zero vector. 
\begin{lemma} (Proposition 2.4.4 in \cite{BLSWZ}) \label{arrange}
There is a homeomorphism $\phi:G(k,\Rn) \to \Arr(k,n)/GL_k$ with $\phi(V) =
[(\alpha \circ \pi_V(e_1), \ldots,
\alpha \circ \pi_V(e_n))]$, where $\pi_V : \Rn \to V$ is orthogonal projection and $\alpha : V \to \Rk$ is any isomorphism.
\end{lemma}

To aid in visualization, we recall in detail the standard map
$m:\Arr(k,n)\to\MacP(k,n)$. Given an arrangement $(v_1,v_2,\ldots,v_n)$, 
the non-zero covectors of the associated oriented matroid are given as follows.
Any oriented 2-dimensional subspace $L$ of ${\mathbb R}^3$ determines a covector  $i
\mapsto +, -, \text{or } 0$ depending on whether $v_i$ is above, is below, or 
 is on $L$. This map $m$ is invariant under the action of $GL_k$ on $\Arr(k,n)$, and
$m[\phi(V)]$ is precisely the oriented matroid associated 
 to each point in $G(k,\Rn)$ mentioned earlier.

For a rank $k$ oriented matroid $M$ on $\{1,2, \dots, n\}$, let $U_M \subset G(k, \Rn)$ be the associated stratum.

\begin{lemma} \label{app:realize} \hfill
\begin{enumerate}
\item For  $M \in \MacP(k,n)$, $\phi (U_M) = m^{-1}(M)/GL_k$.  
\item For  $M \in \MacP(k,n)$, $\phi(\overline{U_M}) = \overline{m^{-1}(M)}/GL_k$.
\end{enumerate}
\end{lemma}

The first statement is clear from the definitions. The second follows from
the first and the fact that $\Arr(k,n) \to \Arr(k,n)/GL_k$ is a principal bundle and hence a closed map.

We will construct a particular oriented matroid $M\in\MacP(3,8)$ and an 
infinite family of oriented matroids $M_i$ whose
properties force any triangulation
of $G(3,\Rinfty)$ refining the matroid stratification
 to have infinitely many simplices in $ U_M$. 

To further help our visualization, consider the set of affine arrangements 
$\AArr(2,n)$, i.e., the set of n-tuples of points in the plane, not all colinear. 
Identify $((x_1,y_1), \dots , (x_n,y_n))$ with $((x_1,y_1,1), \dots,(x_n,y_n,1))$  and thus consider $\AArr(2,n) \subset
\Arr(3,n)$.  Given an affine arrangement $(v_1, \dots v_n)$, the nonzero covectors of the associated oriented matroid are
obtained by considering oriented lines in the affine plane and the positions of the $v_i$ with respect to these lines. For any
$v\in{\mathbb R}^3$ with positive $z$-coordinate $v_z$, let $v'= v_z^{-1} v$. Then the oriented matroid associated to an
arrangement $(v_1,v_2,
\ldots,v_n)$ of
vectors with positive $z$-coordinate is identical to the oriented matroid associated to the affine arrangement $(v_1',v_2',\ldots,v_n')$.
The affine arrangement has the advantage that colinearity, convexity, and intersection properties are determined
by the oriented matroid.

It will be convenient to write elements of $\Nn$ as
$$\{\alpha,\beta,\gamma,\omega,\nu, 
a\}\cup\{b_1,b_2,\ldots\}\cup\{c_1,c_2,\ldots\}
\cup\{d_1,d_2,\dots\}
$$

\begin{figure}[p]
\centerline{\input{a1.pstex_t}}
\caption{The arrangement $A_0$\label{two}}
\end{figure}
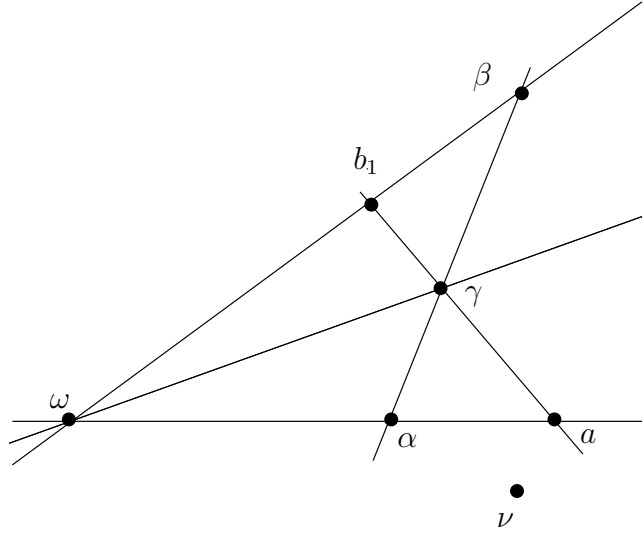

We will define inductively a sequence $A_0 \subset A_1 \subset A_2 \subset \cdots$ of affine 
arrangements.
Consider the affine arrangement $A_0$ pictured
in Figure~\ref{two}.

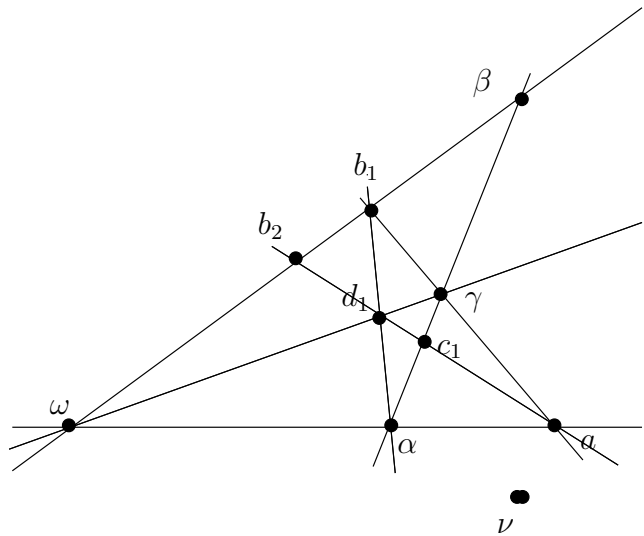
\begin{figure}[p]
\centerline{\input{a2.pstex_t}}
\caption{The arrangement $A_1$\label{three}}
\end{figure}

\begin{figure}[p]
\centerline{\input{a3.pstex_t}}
\caption{The arrangement $A_2$\label{four}}
\end{figure}
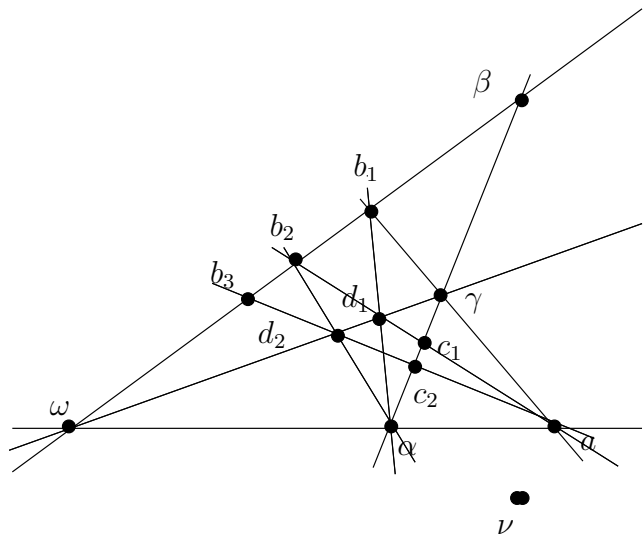

Given an arrangement $A_{n-1}$ with elements $\{\alpha,\beta,\gamma,\omega, \nu,
a\}
\cup\{b_1,b_2,\ldots, b_n\}\cup\{c_1,c_2,\ldots, c_{n-1}\}
\cup\{d_1,d_2,\dots,d_{n-1}\}$, define $A_{n}$ by adding points
$\{d_{n},b_{n+1},c_{n}\}$ to $A_n$ as follows:
\begin{enumerate}
\item Add $d_{n}$ at the intersection of the lines 
$\overline{\omega\gamma}$ and $\overline{\alpha b_{n}}$.
\item Add $b_{n+1}$  at the intersection of the lines 
$\overline{\omega\beta}$ and $\overline{a d_{n}}$.
\item Add $c_{n}$ at the intersection of the lines $\overline{\alpha\beta}$
and $\overline{a b_{n+1}}$.
\end{enumerate}
For instance,  $A_1$ is pictured in Figure~\ref{three} and $A_2$ is pictured 
in Figure~\ref{four}.

An induction on $n$ shows that the $c_i$ are all distinct in each $A_n$.

For each positive integer $i$, let $M_i$ be the oriented matroid 
associated to the arrangement obtained from $A_{i}$ by 
changing the name of $c_i$ to $\delta$.
Note that for any realization of $M_i$ in ${\mathbb R}^3$, the
corresponding realization in affine space is determined by
the positions of $\{\alpha,  a, b_1,\beta, \nu\}.$

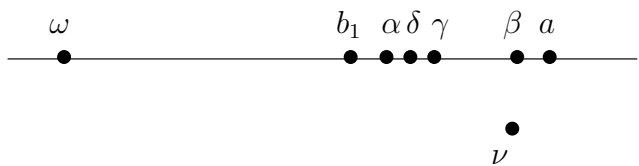
\begin{figure}[p]
\input{m.pstex_t}
\caption{A realization of $M$\label{one}}
\end{figure}

Finally, let $M$ be the oriented matroid associated to the affine
arrangement shown
in Figure~\ref{one}. 

The following lemma says that
the stratum $U_M$ intersects the closure of each $U_{M_i}$, and these
intersections are disjoint.

\begin{lemma}\label{strata} \hfill
\begin{enumerate}
\item $\overline{m^{-1}(M_i)}\cap m^{-1}(M)
\neq\emptyset$ for every $i > 0$, and
\item $\overline{m^{-1}(M_i)}\cap\overline{m^{-1}(M_j)}\cap m^{-1}(M)=\emptyset$
for every $i\neq j$ where $i,j>0$.
\end{enumerate}
\end{lemma}

\begin{proof}
For the first statement, we need a convergent sequence $A_i^1, A_i^2, \dots$ of elements of $\Arr(3,\infty)$, each of which
represents the oriented matroid $M_i$, and whose limit is in $\Arr(3,8)$ and represents $M$.   This sequence is defined by
closing up the angle $\angle \beta \omega a$, leaving the points 
 $\{\alpha,\beta, \gamma, \delta, \omega, \nu, a, b_1\}$ all at height 1 in
$\R^3$ and in the right order in the limit.  Meanwhile, the realizations 
of the remaining points are obtained by letting each $z$-coordinate be
$1/n$ while maintaining the colinearity and intersection properties 
determined by $M_i$.

The second statement is proven using two facts from elementary projective
geometry:
\begin{enumerate}
\item For every two $(n+1)$-tuples $\{x_0,\ldots,x_n\}$ and $\{y_0,\ldots,y_n\}$
of points in general position in affine $n$-dimensional space, there exists a
unique affine automorphism\footnote{An affine automorphism is the composite of a linear automorphism with a
translation.  These are the bijections of affine space which take lines to lines.} taking each
$x_i$ to $y_i$.
\item If $(a,b,c,d)$ are four points on a line in the affine plane
(in the order given), their
\textbf{cross-ratio}
$$  \frac{\|a-c\|}{\|b-c\|}\frac{\|b-d\|}{\|a-d\|}$$
is invariant under affine automorphism.
\end{enumerate}

Now assume by way of contradiction there exists an affine arrangement $B\in
\overline{m^{-1}(M_i)}\cap\overline{m^{-1}(M_j)}\cap m^{-1}(M)$ for some
$i\neq j$. We will compute the cross-ratio of the points $(\alpha,\delta,\gamma,\beta)$ in $B$ in two different ways (an
$i$-way and a $j$-way) and come up with a contradiction.  

Construct a sequence
$B_i^1, B_i^2,\ldots$ of arrangements in 
$m^{-1}(M_i)$ such that
\begin{itemize}
\item the sequence converges to $B$,
\item the elements $\{\alpha,\beta, \gamma,\delta, \omega, \nu, a, b_1\}$ all
have $z$-coordinate 1 in each $B_i^n$, and 
\item the subarrangements $(\alpha,\beta, \gamma,\delta, a, b_1)$ in each
$B_i^n$ are all projectively equivalent. (That is, for every $n_1, n_2$ there 
is an  affine automorphism of the plane taking the 
subarrangement in $B_i^{n_1}$ to the subarrangement in $B_i^{n_2}$.)
\end{itemize}
Fix a small $\epsilon$.
Define $B_i^n$  to be the unique realization of $M_i$ with
\begin{enumerate}
\item $\alpha$, $a$, $\omega$, and $\nu$ in the same positions as in $B$,
\item $\beta$ the point at distance
$\epsilon/n$ above the position of $\beta$ in $B$,
\item $b_1$ determined by requiring the 1-dimensional affine arrangement
$\{\omega, b_1,\beta\}$ in $B_i^n$ to be projectively equivalent to the
corresponding arrangement in $B$, and
\item All other points determined by the colinearity and convexity 
conditions of $M_i$, together with the condition that the $z$-coordinates
of $\gamma$ and $\delta$ are 1 and the $z$-coordinates of all remaining 
points are $1/n$.
\end{enumerate}
That all elements of the arrangements $B_i^n$ except 
$\gamma$ and $\delta$ converge to the
corresponding elements of $B$ is clear. We get convergence of 
$\gamma$ and $\delta$ by noting that there exists some sequence $C_i^1,
C_i^2,\ldots$ in
$m^{-1}(M_i)$ converging to $B$. As $n$ increases, the elements $\{\alpha,\beta,
 \omega, \nu, a, b_1\}$ in $B_i^n$ converge to the corresponding elements of
$C_i^n$. Since the positions of $\gamma$ and $\delta$ are determined by the
positions of $\{\alpha,\beta, a, b_1\}$, $\gamma$ and  $\delta$ converge as well.

Note that the subarrangements  $\{\alpha,\beta, \gamma,\delta, a,
b_1\}$in the $B_i^n$ are all projectively equivalent (by the affine 
automorphism fixing $\omega$ and $a$ and mapping the corresponding $\beta$ to
each other). Thus the cross-ratio $\cross(i)$ of $(\alpha,\delta,\gamma,\beta)$ is
the same in all the $B_i^n$, and so $\cross(i)$
is the cross-ratio 
 of $(\alpha,\delta,\gamma,\beta)$ in $B$. 

Similarly, we get a sequence $B_j^1,B_j^2,\ldots$ in $m^{-1}(M_j)$
and calculate the cross-ratio of $(\alpha,\delta,\gamma,\beta)$ in $B$
to be $\cross(j)$. Thus $\cross(i)=\cross(j)$.
On the other hand, consider the affine automorphism of the plane 
fixing the points $\omega$ and $a$ in $B$ and taking the point $\beta$ in $B_i^1$ to the point $\beta$ in $B_j^1$. This sends the subarrangement
$(\alpha,\beta,a,b_1)$ in $B_i^1$ to the corresponding subarrangement of 
$B_j^1$, hence sends the point $\gamma$ in $B_i^1$ to the point $\gamma$ in 
$B_j^1$. But it does \emph{not} send the point $\delta$  in $B_i^1$ to the point  $\delta$  in $B_j^1$, since $c_i\neq c_j$ in
$A_{\max\{i,j\}}$.  Hence $\cross(i)\neq\cross(j)$, a contradiction.
\end{proof}

\begin{proof}[Proof of Main Theorem] We assume that there is a tame triangulation of $G(3,\R^\infty)$ and reach a
contradiction.  By Lemma \ref{strata} and Lemma \ref{app:realize}, $U_M ~\cap \overline{U_{M_i}} \neq \emptyset$ and 
$U_M ~\cap \overline{U_{M_i}} ~\cap \overline{U_{M_j}}=\emptyset$ for every $i$ and $j$.  Choose a
sequence of simplices $\sigma_1, \sigma_2,\dots$ so that 
$$\text{int } \sigma_i ~\cap~ U_M ~\cap \overline{U_{M_i}} \neq \emptyset.
$$
Then there exists a sequence of simplices $\tau_1, \tau_2,\dots$ so that $\text{ int }\tau_i \subset
U_{M_i}$ and $\sigma_i$ is a face of $\tau_i$.  Then $\text{ int } \sigma_i \subset U_M ~\cap
\overline{U_{M_i}} $ so, by part 2 of the above lemma, the $\sigma_i$'s are distinct.  Thus there are an
infinite number of simplices in a compact set $\overline{U_M} \subset G(3,\R^8)$, which is a contradiction.
\end{proof}

\section{Generalizations}\label{generalize}

Our main theorem can be generalized in two different ways: generalize to partitions more general than a triangulation, and to
stratifications more general than the matroid stratification.

\begin{defn} A \textbf{weak subdivision} of a partition $P=\{U_i: i\in I\}$ of a space $X$ is a partition
$Q$ of $X$ such that
\begin{itemize} 
\item $Q$ \emph{refines} $P$, i.e.,  every $U$ in $P$ is the union of elements of $Q$,
\item Each element of $Q$ is \emph{connected},
\item $Q$ is \emph{locally finite}, i.e. every compact set $K$ of $X$ intersects only a finite number of elements of $Q$,
\item $Q$ is \emph{normal}, i.e., if $U$ and $V$ are elements of $Q$ and
$U\cap\overline{V}\neq\emptyset$ then $U\subseteq \overline{V}$.
\end{itemize}
\end{defn}
If $Q$ consists of the interiors of simplices in a triangulation of $X$, and $Q$ refines a partition $P$, then $Q$ is a weak 
subdivision of
$P$.  However, a CW decomposition of $X$ refining $P$ need not be a weak subdivision; normality may not hold.

If $M$ and $M'$ are (oriented) matroids on a set $E$ and if every covector of $M'$ is a covector of $M$, then one says
that there is a \textbf{strong map $M \to M'$}.  Define the \textbf{combinatorial Grassmannian}
$$
\Gamma(k,M) =\{M': \text{rank }M' = k \text{ and } M  \to M' \}.
$$
If $M$ is the (oriented) matroid associated to a 
collection of vectors $\{v_1, v_2,
\dots, v_m\}$ spanning
$\Rn$, define
$$
\mu_M : G(k,\Rn) \to \Gamma(k,M)
$$
by sending $V$ to the (oriented) matroid 
with covectors
$$\{i \mapsto \text{sign}\langle v, v_i \rangle\}_{v\in V}.
$$
The
point inverse images are called the \textbf{generalized (oriented) matroid strata}.  This stratification comes up in the
study of extension spaces of oriented matroids (cf.~\cite{SZ}). 

To get a stratification of the infinite Grassmannian, one needs some compatibility between the matroids.  Let 
$ A_1 \subseteq A_2 \subseteq A_3 \subseteq \dots \subset \Rinfty$ be a sequence of finite sets, whose union spans
$\Rinfty$, and so that every element of $A_{i+1} - A_i$ has all of its first $i$ coordinates zero.  Then the associated
oriented matroids $M_1, M_2, \cdots$ 
satisfy
\begin{itemize}
\item $M_{i+1}\rightarrow M_i$ for all $i$, and so there are inclusions
$\Gamma(k, M_i)\to\Gamma(k, M_{i+1})$, and
\item the associated maps $\mu_{M_i} : G(k, \R^{n_i}) \to \Gamma(k,M_i)$
commute with the inclusions of real resp. combinatorial Grassmannians.
\end{itemize}
Thus the maps $\mu_{M_i}$ give a generalized (oriented) matroid stratification
of $G(k,\Rinfty)$.  By choosing $B_1 \subseteq B_2 \subseteq B_3 \subseteq \dots$ so that the $B_i \subseteq A_i$ and the
$B_i$ are linearly independent, we see that the generalized (oriented) matroid stratification is essentially a refinement
of the oriented matroid stratification.  Thus the our proof of our main theorem showed the following result.

\begin{theorem}  There is no weak subdivision of any generalized (oriented) matroid stratification of $G(k,\Rinfty)$ for
$k \geq 3$.
\end{theorem}

\section{Tame Triangulations and matroid bundles}\label{context}

The point of this section is to give a context for tame triangulations, to show how they allow the process of
\emph{combinatorialization}, the passage from topological structures to combinatorial ones.  This occurs in two related ways, in constructing
maps from real to combinatorial Grassmannians, and in passing from vector bundles to matroid bundles.  For more on this
see \cite{Z2} and  \cite{bundles}.  
Our main theorem forces delicate constructions in going from
the finite to the infinite dimensional case
in \cite{Z2}.

Let $\pi : E \to B$ be a rank $k$ vector bundle over a simplicial complex.  Assume the fibers $F_b = \pi^{-1}(b)$ are
equipped with a continuously varying inner product.  If $B$ is finite-dimensional, there is a set of spanning sections $S =
\{s_1, s_2, \dots , s_n\}$.  Then we have a map
$$
\Mm : B \to \MacP(k,n)
$$
sending $b \in B$ to the oriented matroid associated to $\{s_1(b), s_2(b), \dots , s_n(b)\} \subset F_b$. (The reader is
strongly urged to work out the case of the open M\"obius strip mapping to the circle.)  The set of sections $S$ is said to
be \textbf{tame} when $\Mm$ is constant on the interior of simplices, in which case we have a true combinatorial gadget,
a \emph{matroid bundle}.  The existence of a tame triangulation of $G(k,\Rn)$ shows that every rank $k$ vector bundle over an
$n$-dimensional complex has a tame set of  sections after subdividing the base.  This is accomplished by applying the
simplicial approximation theorem to a classifying map $B \to G(k,\Rn)$, and pulling back the canonical sections.

One can think of the MacPhersonian $\MacP(k,n)$ of rank $k$ oriented matroids on $\{1,2, \dots ,n\}$ as a classifying
space for matroid bundles.  It has a partial order given by
$M_1
\geq M_2$ if there is a weak map from $M_1$ to $M_2$.  If $U_{M_1} \cup \overline{U_{M_2}} \neq \emptyset$, then $M_1 \geq
M_2$.  Let $\|\MacP(k,n)\|$ be the geometric realization (= order complex) of this poset.  Let $\mu : G(k,\Rn) \to
\MacP(k,n)$ be the realization map.  Given a tame triangulation of $G(k,\Rn)$, then one can construct a simplicial map
$\tilde{\mu} : G(k,\Rn) \to \|\MacP(k,n)\|$ from the barycentric subdivision of the tame triangulation agreeing with $\mu$
on the vertices. The main result of \cite{Z2} shows that $\tilde{\mu}$ carries the Stiefel-Whitney classes, and hence
 Stiefel-Whitney classes can be defined purely combinatorially.

\vskip 24pt
\noindent\textbf{Acknoledgements:} We thank G\"unter Ziegler for helpful comments.
 
\bibliographystyle{alpha}
\bibliography{biblio}

\end{document}

%% file: a1.pstex_t
\begin{picture}(0,0)%
\epsfig{file=a1.pstex}%
\end{picture}%
\setlength{\unitlength}{3947sp}%
\begingroup\makeatletter\ifx\SetFigFont\undefined
% extract first six characters in \fmtname
\def\x#1#2#3#4#5#6#7\relax{\def\x{#1#2#3#4#5#6}}%
\expandafter\x\fmtname xxxxxx\relax \def\y{splain}%
\ifx\x\y   % LaTeX or SliTeX?
\gdef\SetFigFont#1#2#3{%
  \ifnum #1<17\tiny\else \ifnum #1<20\small\else
  \ifnum #1<24\normalsize\else \ifnum #1<29\large\else
  \ifnum #1<34\Large\else \ifnum #1<41\LARGE\else
     \huge\fi\fi\fi\fi\fi\fi
  \csname #3\endcsname}%
\else
\gdef\SetFigFont#1#2#3{\begingroup
  \count@#1\relax \ifnum 25<\count@\count@25\fi
  \def\x{\endgroup\@setsize\SetFigFont{#2pt}}%
  \expandafter\x
    \csname \romannumeral\the\count@ pt\expandafter\endcsname
    \csname @\romannumeral\the\count@ pt\endcsname
  \csname #3\endcsname}%
\fi
\fi\endgroup
\begin{picture}(3999,3312)(5614,-6061)
\put(8077,-5562){\makebox(0,0)[lb]{\smash{\SetFigFont{12}{14.4}{rm}$\alpha$}}}
\put(5893,-5315){\makebox(0,0)[lb]{\smash{\SetFigFont{12}{14.4}{rm}$\omega$}}}
\put(9226,-5536){\makebox(0,0)[lb]{\smash{\SetFigFont{12}{14.4}{rm}$a$}}}
\put(8551,-3286){\makebox(0,0)[lb]{\smash{\SetFigFont{12}{14.4}{rm}$\beta$}}}
\put(7801,-3811){\makebox(0,0)[lb]{\smash{\SetFigFont{12}{14.4}{rm}$b_1$}}}
\put(8501,-4636){\makebox(0,0)[lb]{\smash{\SetFigFont{12}{14.4}{rm}$\gamma$}}}
\put(8701,-6061){\makebox(0,0)[lb]{\smash{\SetFigFont{12}{14.4}{rm}$\nu$}}}
\end{picture}

%% file: a2.pstex_t
\begin{picture}(0,0)%
\epsfig{file=a2.pstex}%
\end{picture}%
\setlength{\unitlength}{3947sp}%
\begingroup\makeatletter\ifx\SetFigFont\undefined
% extract first six characters in \fmtname
\def\x#1#2#3#4#5#6#7\relax{\def\x{#1#2#3#4#5#6}}%
\expandafter\x\fmtname xxxxxx\relax \def\y{splain}%
\ifx\x\y   % LaTeX or SliTeX?
\gdef\SetFigFont#1#2#3{%
  \ifnum #1<17\tiny\else \ifnum #1<20\small\else
  \ifnum #1<24\normalsize\else \ifnum #1<29\large\else
  \ifnum #1<34\Large\else \ifnum #1<41\LARGE\else
     \huge\fi\fi\fi\fi\fi\fi
  \csname #3\endcsname}%
\else
\gdef\SetFigFont#1#2#3{\begingroup
  \count@#1\relax \ifnum 25<\count@\count@25\fi
  \def\x{\endgroup\@setsize\SetFigFont{#2pt}}%
  \expandafter\x
    \csname \romannumeral\the\count@ pt\expandafter\endcsname
    \csname @\romannumeral\the\count@ pt\endcsname
  \csname #3\endcsname}%
\fi
\fi\endgroup
\begin{picture}(3999,3312)(5614,-6061)
\put(8077,-5562){\makebox(0,0)[lb]{\smash{\SetFigFont{12}{14.4}{rm}$\alpha$}}}
\put(5893,-5315){\makebox(0,0)[lb]{\smash{\SetFigFont{12}{14.4}{rm}$\omega$}}}
\put(9226,-5536){\makebox(0,0)[lb]{\smash{\SetFigFont{12}{14.4}{rm}$a$}}}
\put(8551,-3286){\makebox(0,0)[lb]{\smash{\SetFigFont{12}{14.4}{rm}$\beta$}}}
\put(7801,-3811){\makebox(0,0)[lb]{\smash{\SetFigFont{12}{14.4}{rm}$b_1$}}}
\put(8501,-4636){\makebox(0,0)[lb]{\smash{\SetFigFont{12}{14.4}{rm}$\gamma$}}}
\put(8701,-6061){\makebox(0,0)[lb]{\smash{\SetFigFont{12}{14.4}{rm}$\nu$}}}
\put(7726,-4636){\makebox(0,0)[lb]{\smash{\SetFigFont{12}{14.4}{rm}$d_1$}}}
\put(8326,-4936){\makebox(0,0)[lb]{\smash{\SetFigFont{12}{14.4}{rm}$c_1$}}}
\put(7201,-4186){\makebox(0,0)[lb]{\smash{\SetFigFont{12}{14.4}{rm}$b_2$}}}
\end{picture}

%% file: a3.pstex_t
\begin{picture}(0,0)%
\epsfig{file=a3.pstex}%
\end{picture}%
\setlength{\unitlength}{3947sp}%
\begingroup\makeatletter\ifx\SetFigFont\undefined
% extract first six characters in \fmtname
\def\x#1#2#3#4#5#6#7\relax{\def\x{#1#2#3#4#5#6}}%
\expandafter\x\fmtname xxxxxx\relax \def\y{splain}%
\ifx\x\y   % LaTeX or SliTeX?
\gdef\SetFigFont#1#2#3{%
  \ifnum #1<17\tiny\else \ifnum #1<20\small\else
  \ifnum #1<24\normalsize\else \ifnum #1<29\large\else
  \ifnum #1<34\Large\else \ifnum #1<41\LARGE\else
     \huge\fi\fi\fi\fi\fi\fi
  \csname #3\endcsname}%
\else
\gdef\SetFigFont#1#2#3{\begingroup
  \count@#1\relax \ifnum 25<\count@\count@25\fi
  \def\x{\endgroup\@setsize\SetFigFont{#2pt}}%
  \expandafter\x
    \csname \romannumeral\the\count@ pt\expandafter\endcsname
    \csname @\romannumeral\the\count@ pt\endcsname
  \csname #3\endcsname}%
\fi
\fi\endgroup
\begin{picture}(3999,3312)(5614,-6061)
\put(8077,-5562){\makebox(0,0)[lb]{\smash{\SetFigFont{12}{14.4}{rm}$\alpha$}}}
\put(5893,-5315){\makebox(0,0)[lb]{\smash{\SetFigFont{12}{14.4}{rm}$\omega$}}}
\put(9226,-5536){\makebox(0,0)[lb]{\smash{\SetFigFont{12}{14.4}{rm}$a$}}}
\put(8551,-3286){\makebox(0,0)[lb]{\smash{\SetFigFont{12}{14.4}{rm}$\beta$}}}
\put(7801,-3811){\makebox(0,0)[lb]{\smash{\SetFigFont{12}{14.4}{rm}$b_1$}}}
\put(8501,-4636){\makebox(0,0)[lb]{\smash{\SetFigFont{12}{14.4}{rm}$\gamma$}}}
\put(8701,-6061){\makebox(0,0)[lb]{\smash{\SetFigFont{12}{14.4}{rm}$\nu$}}}
\put(7726,-4636){\makebox(0,0)[lb]{\smash{\SetFigFont{12}{14.4}{rm}$d_1$}}}
\put(8326,-4936){\makebox(0,0)[lb]{\smash{\SetFigFont{12}{14.4}{rm}$c_1$}}}
\put(6901,-4486){\makebox(0,0)[lb]{\smash{\SetFigFont{12}{14.4}{rm}$b_3$}}}
\put(7201,-4861){\makebox(0,0)[lb]{\smash{\SetFigFont{12}{14.4}{rm}$d_2$}}}
\put(8176,-5236){\makebox(0,0)[lb]{\smash{\SetFigFont{12}{14.4}{rm}$c_2$}}}
\put(7276,-4186){\makebox(0,0)[lb]{\smash{\SetFigFont{12}{14.4}{rm}$b_2$}}}
\end{picture}

%% file: m.pstex_t
\begin{picture}(0,0)%
\psfig{file=m.pstex}%
\end{picture}%
\setlength{\unitlength}{3947sp}%
\begingroup\makeatletter\ifx\SetFigFont\undefined
% extract first six characters in \fmtname
\def\x#1#2#3#4#5#6#7\relax{\def\x{#1#2#3#4#5#6}}%
\expandafter\x\fmtname xxxxxx\relax \def\y{splain}%
\ifx\x\y   % LaTeX or SliTeX?
\gdef\SetFigFont#1#2#3{%
  \ifnum #1<17\tiny\else \ifnum #1<20\small\else
  \ifnum #1<24\normalsize\else \ifnum #1<29\large\else
  \ifnum #1<34\Large\else \ifnum #1<41\LARGE\else
     \huge\fi\fi\fi\fi\fi\fi
  \csname #3\endcsname}%
\else
\gdef\SetFigFont#1#2#3{\begingroup
  \count@#1\relax \ifnum 25<\count@\count@25\fi
  \def\x{\endgroup\@setsize\SetFigFont{#2pt}}%
  \expandafter\x
    \csname \romannumeral\the\count@ pt\expandafter\endcsname
    \csname @\romannumeral\the\count@ pt\endcsname
  \csname #3\endcsname}%
\fi
\fi\endgroup
\begin{picture}(3979,1020)(5634,-6061)
\put(8701,-6061){\makebox(0,0)[lb]{\smash{\SetFigFont{12}{14.4}{rm}$\nu$}}}
\put(9001,-5236){\makebox(0,0)[lb]{\smash{\SetFigFont{12}{14.4}{rm}$a$}}}
\put(5926,-5236){\makebox(0,0)[lb]{\smash{\SetFigFont{12}{14.4}{rm}$\omega$}}}
\put(8326,-5236){\makebox(0,0)[lb]{\smash{\SetFigFont{12}{14.4}{rm}$\gamma$}}}
\put(8776,-5236){\makebox(0,0)[lb]{\smash{\SetFigFont{12}{14.4}{rm}$\beta$}}}
\put(7726,-5236){\makebox(0,0)[lb]{\smash{\SetFigFont{12}{14.4}{rm}$b_1$}}}
\put(8011,-5236){\makebox(0,0)[lb]{\smash{\SetFigFont{12}{14.4}{rm}$\alpha$}}}
\put(8166,-5236){\makebox(0,0)[lb]{\smash{\SetFigFont{12}{14.4}{rm}$\delta$}}}
\end{picture}